\documentclass[12pt]{amsart}
\usepackage{amssymb,latexsym}
\usepackage{enumerate}

\makeatletter
\@namedef{subjclassname@2010}{%
  \textup{2010} Mathematics Subject Classification}
\makeatother

\newtheorem{thm}{Theorem}[section]
\newtheorem{cor}[thm]{Corollary}
\newtheorem{lem}[thm]{Lemma}
\newtheorem{prob}[thm]{Problem}
\newtheorem{prop}[thm]{Proposition}

\theoremstyle{definition}
\newtheorem{defin}[thm]{Definition}
\newtheorem{rem}[thm]{Remark}

\numberwithin{equation}{section}

\begin{document}


\baselineskip=17pt



\author[A. Elkhadiri]{Abdelhafed Elkhadiri}
\address{Department of Mathematics\\ Faculty   of Sciences\\
University Ibn Tofail B.P. 133, K\'enitra, Morocco}
\email{kabdelhafed@hotmail.com}
\thanks{ This work was partially supported by PARS MI33}

\date{}

\subjclass[2010]{Primary 14B25, 14B10, 32B05, 32B20; Secondary
03C10}

\keywords{ Weierstrass division theorem, semi-analytic sets,
o-minimal structures  }

\title[On  connected components..] {On  connected components of some globally semi-analytic
sets}

\begin{abstract}
We isolate a class, say $\mathcal{A}$, of global real analytic
functions such that, each global semi-analytic set defined by
$\mathcal{A}$ has only  finitely   many  connected components and
each component is also a global semi-analytic set defined by
$\mathcal{A}$.
\end{abstract}
\maketitle

\section{Introduction}
A subset of $\mathbb{R}^n$ is called semi-algebraic if it can be
represented as a (finite ) boolean combination of sets of the form
$\{ x\in\mathbb{R}^n\,\mid P(x)=0\}$, $\{ x\in\mathbb{R}^n\,\mid
Q(x)\neq 0\}$, where $P(x)$ and $Q(x)$ are $n$ variables polynomials
with real coefficients. A fundamental result of Tarski-Seidenberg
says that the projection of a semi-algebraic set is a semi-algebraic
set. This result is known to logicians as quantifiers elimination
for the ordered ring structure on $\mathbb{R}$, say
$\overline{\mathbb{R}}$. Immediate consequence are the facts that
the closure, interior and boundary of a semi-algebraic set are
semi-algebraic. The result of Tarski-Seidenberg is also the basis
for many inductive arguments in semi-algebraic geometry where a
desired property of a given semi-algebraic set is inferred from the
same property of projection of the set into lower dimension. See for
example \cite{Hiro} and \cite{Coste} for triangulation of bounded
semi-algebraic sets.\\ An other important property of semi-algebraic
sets is that they have only a finite number of  connected
components, and each of them  is also  semi-algebraic. This follows
from Collins \cite{Colline}, in which a new decision method for the
reals is constructed, much more time efficient than Tarski's. Here
we are interested in Collins' key {\it{geometric}} idea, which he
calls "{ \it{cylindric decomposition}}". This idea is used by
Bierston and Milman \cite{Bier} to establish some  results of
Lojasiewicz \cite{Loja} on the germs of
semi-analytic sets.\\
Let $\mathcal{R}$ be an  arbitrary but fixed o-minimal expansion of
the ordered field $\mathbb{R}$. For more details about o-minimal structures over the field of reals,
we refer the reader to \cite{LVD}.
Consider $\Omega\subset\mathbb{R}^n$
an open definable set. In this paper we isolate a subalgebra, say
$\mathcal{A}(\Omega)$, of the algebra of real analytic definable
functions on $\Omega$, and we define our global semi-analytic sets
as follow:  a  set $ A\subset \Omega$ is called $\mathcal{A}$- semi
analytic set if there exists a finite number of open definable sets
$\Omega_j\subset \Omega,\,1\leq j\leq s,$ such that:
\begin{enumerate}
\item  $\Omega=\bigsqcup_{j=1}^s\Omega_j$
\item   $ A\cap\Omega_j$ is a finite union of sets of the form:
 \[\{ x\in\Omega /\, \varphi_0(x)=0,\,\varphi_1(x)>0,\ldots,\varphi_k
>0\},\]
\end{enumerate}

with $\varphi_0,\varphi_1,\ldots,\varphi_k\in\mathcal{A}(\Omega_j)$.

The main result of this paper is a semi-global Weierstrass  Preparation  Theorem
 in such algebra $\mathcal{A}(\Omega)$, see corollary 5.6. As an application, we prove that every $\mathcal{A}$- semi analytic set has a finite number of connected components, and each of them is also  $\mathcal{A}$- semi analytic.
\section{Topological Noetherianity  }
$\mathcal{R}=(\mathbb{R}, <, +,.,-,0,1,\ldots)$ will, in this paper,
denote an arbitrary, but fixed, o-minimal expansion of the ordered
field $\mathbb{R}$. Definable means definable in $\mathcal{R}$ with
parameters from $\mathbb{R}$. Let $\Omega\subset\mathbb{R}^n$ be an
open set, we denote by $ \mathcal{H}(\Omega)$ the algebra of real
analytic functions on $\Omega$. If $\Omega$ is definable, we denote
by $\mathcal {O}(\Omega)\subset \mathcal{H}(\Omega)$  the set of all
$f\in \mathcal{H}(\Omega)$ definable. Clearly, $\mathcal
{O}(\Omega)$ is an $\mathbb{R}$-subalgebra closed under derivative
and $\mathbb{R}[x_1,\ldots,x_n]_{|\Omega}\subset\mathcal
{O}(\Omega)$, where  $\mathbb{R}[x_1,\ldots,x_n]$ is the ring of
polynomials. If $SM \mathcal {O}(\Omega)$ is the maximal spectrum of
$\mathcal {O}(\Omega)$, we have then an injection
\[\Omega\rightarrow SM \mathcal {O}(\Omega),\] every $x\in \Omega $
is identified with the maximal ideal: \[\underline{m}_x=\{ f\in
\mathcal {O}(\Omega)\,/\, f(x)=0\}.\] We denote by
$\Omega(\mathcal{O})$ the topological space $\Omega$ with the
induced topology of $SM \mathcal {O}(\Omega)$.
\begin{lem}
Let $\Omega\subset \mathbb{R}$ be an open definable set, then the
ring $\mathcal {O}(\Omega)$ is Noetherian.
\end{lem}
\underline{Proof}. Since $\Omega$ is a finite union of open disjoint
intervals: $\Omega_1,\ldots,\Omega_q$, we have
$\mathcal{O}(\Omega)=\mathcal{O}(\Omega_1)\oplus\ldots\oplus
\mathcal{O}(\Omega_q).$ We can then suppose, for the proof, that
$\Omega$ is an open interval. Let $F\in\mathcal{O}(\Omega)-\{0\}$,
since the number of real zeros of $F$ in $\Omega$ is finite, there
exists a polynomial $Q_F\in\mathbb{R}[X]$ such that $F=Q_F\Psi_F$,
where $\Psi_F\in \mathcal{H}(\Omega)$ and $ \Psi_F(x)\neq
0,\,\,\forall x\in\Omega$. It is clear that
$\Psi_F\in\mathcal{O}(\Omega)$.\\ Now let $
I\subset\mathcal{O}(\Omega)$ be an ideal, we consider the ideal
$J\subset\mathbb{R}[x_1,\ldots,x_n]$ generated by all $Q_F,\, F\in
I$. Let $\{Q_{F_1},\ldots,Q_{F_s}\}$ be a system of generators of
the ideal $J$. It is clear that the system $\{ F_1,\ldots,F_s\}$
generates the ideal $I$.
\begin{rem}
If $\Omega= \mathbb{R}$, the algebra
$\mathbb{R}[x,\,\sin\,x,\,\cos\,x]$ is Noetherian, but the functions
$\mathbb{R}\ni x\rightarrow \sin\,x$ and $\mathbb{R}\ni x\rightarrow \cos\,x$ are not definable
in any o-minimal structure.
\end{rem}

Let $\Omega\subset \mathbb{R}^n$ be an open set,
if $ I\subset \mathcal{H}(\Omega)$ is an a ideal, we put:
\[ V(I)=\{x\in \Omega\,/f(x)=0,\,\forall f\in I\}.\]
If $I$ is generated by one element $h\in \mathcal{H}(\Omega)$, we write $V(h)$ instead of $ V((h)\mathcal{H}(\Omega))$.
Finely, we denote by   $Reg V(h)$
the set of regular points of $V(h)$ i.e. the set of all $x\in V(h)$
such that $ V(h)$ is an embedded submanifold  in a neighborhood of $x$.\\
We consider the following property for the structure $\mathcal{R}$.
\subsection{}
\hspace{2cm}{\it{ For all $f\in \mathcal{O}(\Omega)$, $Reg V(f)$ is definable in $\mathcal{R}$}}.\\
We recall that if the structure $\mathcal{R}$ has the analytic cell decomposition, see \cite{LVD1}, then $\mathcal{R}$ satisfies the property
2.1. We can find examples of such structures in \cite{Miller}. We note that not every 0-minimal structure has analytic cell decomposition see \cite{le gal}.
However, most known 0-minimal expansions of the real numbers admit analytic cell decomposition. Both $\overline{\mathbb{R}}=(\mathbb{R}, <, +,.,-,0,1)$ and $(\overline{\mathbb{R}}, \exp)$, where $x\rightarrow \exp x$ is the global exponential function, admit analytic cell decomposition. More generally the Pfaffian closure of an o-minimal structure preserves analytic cell decomposition.
\begin{rem}
In the case where the structure $\mathcal{R}=\overline{\mathbb{R}}=(\mathbb{R}, <,
+,.,-,0,1)$, if $\Omega\subset \mathbb{R}^n$ is an open definable set,
then $\Omega$ is semi-algebraic. In this setting
$\mathcal{O}(\Omega)$ is the algebra of Nash functions, we know by
\cite{Ris} that $\mathcal{O}(\Omega)$ is Noetherian.
\end{rem}

\begin{prob}
Let $\Omega\subset \mathbb{R}^n,\, n >1$, be an open definable set, we suppose that the structure satisfies the property 2.1.
Is the algebra $\mathcal{O}(\Omega)$ is Noetherian?
\end{prob}

The following theorem is the topological version of Noetherianity of
$\mathcal{O}(\Omega)$.
\begin{thm}
Suppose that the structure $\mathcal{R}$ satisfies the property 2.1, then the topological space $\Omega(\mathcal{O})$ is a Noetherian space, in order
words, any decreasing sequence of closed sets of
$\Omega(\mathcal{O})$ stabilizes.
\end{thm}
\underline{Proof}. Let $F\subset\Omega(\mathcal{O})$ be a closed
set. There exists $ I\subset\mathcal{O}(\Omega)$ an ideal such that
$F=V(I):=\{ x\in\Omega\,/\, f(x)=0,\,\forall f\in I\}$. First, we
will show that $ V(I)$ can be defined by one equation, i.e. there is
$h\in I$ such that \[V(I)= V(h):=\{x\in\Omega\,/\,h(x)=0\}.\] Let
$g\in I-\{0\}$, we have $ V(I)\subset V(g)$.  Since  $Reg V(g)$ is definable, we know that  $Reg V(g)$ has a finite
number of connected components: $\Gamma_1,\ldots,\Gamma_s$. We put
\[\mu_1=max\{ dim\Gamma_j\,/\, \Gamma_j\nsubseteq V(I)\}.\] If
$\Gamma_1,\ldots,\Gamma_\nu$ are the connected components of $Reg
V(g)$ such that \[dim\Gamma_l = \mu_1,\, \Gamma_l\nsubseteq V(I),\,l=1,\ldots,\nu,\] for each
$l=1,\ldots,\nu$, there exists $h_l\in I$ such that
${h_l}_{|\Gamma_l}\neq 0$. We put $h=\sum_{l=1}^\nu h_l^2\in I$. We
have \[ V(I)\subset V(\psi)\varsubsetneq V(g),\] where $\psi=h^2+g^2$.\\
If we put \[\mu_2=max\{ dim\Gamma_k\,/\,\Gamma_k \,\mbox{connected
component of}\,Reg \psi,\,
 \Gamma_k\nsubseteq V(I)\},\] we have $\mu_2 <\mu_1$.\\ By continuing
 this processus with $\psi$ and so on, we see then, there is
 $\varphi\in I$ such that $ V(I)= V(\varphi)$.\\
 Now let $(F_j)_{j\in J}$ be a decreasing sequence of closed sets of the
 topological space $\Omega(\mathcal{O})$. For each $j\in J$, there
 exists $\varphi_j\in\mathcal{O}(\Omega)$ such that $F_j=
 V(\varphi_j)$. For each $V(\varphi_j)$ we associate a $n+1$-tuple
 $\nu_j=(\nu_{j,n},\nu_{j,n-1},\ldots,\nu_{j,0})\in\mathbb{N}^{n+1}$,
 where $\nu_{j,k}$ is the number of connected components of $Reg
 V(\varphi_j)$ of dimension $k$. If we consider the lexicographic order
 on $\mathbb{N}^{n+1}$,  we have $\nu_j < \nu_i$ if $V(\varphi_j)\subsetneqq
 V(\varphi_i)$. Hence the sequence $(F_j)_{j\in J}$ stabilizes.
\begin{prop}
With the same hypothesis as theorem 2.5, the mapping $\Omega\rightarrow SM \mathcal {O}(\Omega)$ is a
bijection.
\end{prop}
\underline{Proof}. Let $\underline{m}\in SM \mathcal {O}(\Omega)$
and suppose that for each $x\in\Omega$, $\underline{m}\nsubseteq
\underline{m}_x$. If $f\in\underline{m}$, there is $x_1\in\Omega$
such that $f(x_1)=0$ ( if not $f$ will be invertible in
$\mathcal{O}(\Omega)$). Since $\underline{m}\nsubseteq
\underline{m}_{x_1}$, there exists $g\in\underline{m}$ with $
g(x_1)\neq 0$. We have then $ V(h)\subsetneqq V(f)$, where $h=
f^2+g^2$. We remark that $ V(h)\neq\emptyset$, since
$h\in\underline{m}$. We pick $x_2\in V(h)$ and by using again
$\underline{m}\nsubseteq \underline{m}_{x_2}$, there exits
$k\in\underline{m}$ with $ k(x_2)\neq 0$. By continuing this
processus, we construct an infinite decreasing sequence of closed
sets of $\Omega(\mathcal{O})$, which is a contradiction. Hence there
exists $x\in\Omega$ such that $\underline{m}\subset \underline{m}_x$
i.e. $\underline{m}= \underline{m}_x$.\\
 In all the following, $\mathcal{R}$ is an o-minimal expansion of the field of real numbers satisfying the property 2.1.
  We define  a subalgebra, $\mathcal{A}(\Omega)\subset
\mathcal{O}(\Omega)$ say, which will be used to define our family of
global semi-analytic sets.
\subsection{$\mathbb{C}$-definable functions}   We
begin by some notations and
conventions. \\
 We identify $\mathbb{C}$ with
$\mathbb{R}^2$, and view subsets of $\mathbb{C}^n$ as subsets of
$\mathbb{R}^{2n}$. Under this identification, functions from
$\mathbb{C}^n$ into $\mathbb{C}$ become maps from $\mathbb{R}^{2n}$
into $\mathbb{R}^2$. The following facts will be constantly used
without further mention. If $\Omega\subset\mathbb{R}^n$ is an open
set and  if $f\in\mathcal{H}(\Omega)$, there exist an open set
$\tilde{\Omega}\subset\mathbb{C}^n$ and $\tilde{f}$ a  holomorphic
function on $\tilde{\Omega}$, such that $\tilde{\Omega}\cap
\mathbb{R}^n= \Omega$ and $\tilde{f}_{|\Omega}=f$. We denote by
$Re\tilde{f}$ and $Im\tilde{f}$ the real and imaginary part of
$\tilde{f}$ respectively.
\begin{defin}
Let $\Omega\subset \mathbb{R}^n$ be an open definable set, and a
function $ f\in\mathcal{O}(\Omega)$. We call $f${\it{
$\mathbb{C}$-definable,}} if the holomorphic function $\tilde{f}$ is
definable i.e. $Re\tilde{f}$ and $Im\tilde{f}$ are definable.
\end{defin}
\begin{rem} 
We note that, not every definable function is $\mathbb{C}$-definable, for example, if $\mathcal{R}= ( \overline{\mathbb{R}}, \exp)$,
the function $x\rightarrow \exp x$ is definable but $Re\tilde{f}(x)= \cos x$ and $Im\tilde{f}(x)= \sin x$ are not definable.
\end{rem}
We denote by $ \mathcal{A}(\Omega)\subset \mathcal{O}(\Omega)$ the
subalgebra of $\mathbb{C}$-definable functions. clearly $
\mathcal{A}(\Omega)$ is closed under derivative and
$\mathbb{R}[x_1,\ldots,x_n]_{|\Omega}\subset\mathcal{A}(\Omega)$.
\section{ Generic Division}
We denote by $z=(z_1,\ldots,z_n)$ and $v=(v_1,\ldots,v_p)$ the
coordinates of $\mathbb{C}^n$ and $\mathbb{C}^p$ respectively. We
put $z'=(z_1,\ldots,z_{n-1})$, $z_j=x_j+iy_j, \, 1\leq j\leq n$,
$v_j=t_j+i\tau_j, \,1\leq j\leq p$, with $i^2=-1$.\\ If $n',p'\in
\mathbb{N},\,n'\leq n,\,p'\leq p$, we denote by $\pi_{(n\times
p,n'\times p')}:\mathbb{C}^n\times \mathbb{C}^p\rightarrow
\mathbb{C}^{n'}\times \mathbb{C}^{p'}$ the
projection defined by
$ \pi_{(n\times p,n'\times p')}(z_1,\ldots,z_n,v_1,\ldots,v_p) = (z_1,\ldots,z_{n'},v_1,\ldots,v_{p'}).$
We consider the mapping $\sigma:\mathbb{C}^p\rightarrow
\mathbb{C}^p$ defined by $ \sigma(v)=( -\sigma_1(v),\ldots,
(-1)^p\sigma_p(v))$, where $\sigma_l$ is the $l$-th elementary
symmetric function of the variables $(v_1,\ldots,v_p)$, $1\leq l\leq
p$. \\ For each permutation $\alpha$ of the set $\{1,\ldots,p\}$, we
associate a $ \mathbb{C}$-linear isomorphism
$\pi_\alpha:\mathbb{C}^n\times \mathbb{C}^p\rightarrow
\mathbb{C}^n\times\mathbb{C}^p$ as follow:
\[\pi_\alpha( z,v_1,\ldots,v_p)= ( z,
v_{\alpha(1)},\ldots,v_{\alpha(p)}) .\]
A set
$\Lambda\subset\mathbb{C}^n\times \mathbb{C}^p$ is called symmetric,
if $\pi_\alpha(\Lambda)= \Lambda$. Let
$\Lambda\subset\mathbb{C}^n\times \mathbb{C}^p$ be a symmetric set.
A fonction $ F:\Lambda\subset\mathbb{C}^n\times
\mathbb{C}^p\rightarrow \mathbb{C}$ is called symmetric, if $F\circ
\pi_\alpha= F$.\\
For each $v=(v_1,\ldots,v_n)\in\mathbb{C}^p$,
$P(v,z_n)=z_n^p+\sum_{j=1}^pv_jz_n^{p-j}$ is called a generic
polynomial in $z_n$ of degree $p$. We put $|v|=max_{j=1}^p|v_j|$.
\begin{lem}\cite{Coste}
Let $\xi$ be a root of the polynomial
$P(v,z_n)=z_n^p+\sum_{j=1}^pv_jz_n^{p-j}$, then $|\xi|\leq 2\,
max_{j=1}^p|v_j|^{\frac{1}{j}}$.
\end{lem}
Let $\Omega\subset\mathbb{R}^n$ be an open definable set and
$g\in\mathcal{A}(\Omega)$, there exist an open definable
$\tilde{\Omega}\subset\mathbb{C}^n$ and a definable holomorphic
function $ \tilde{g}$ on $\tilde{\Omega}$ such that
$\tilde{g}_{|\Omega}=g$ and $\tilde{\Omega}\cap\mathbb{R}^n=\Omega$.
We consider:
\[ \tilde{W}:=\{ (z',z_n,v)\in \tilde{\Omega}\times \mathbb{C}^p/\,
\forall \xi\, \mbox{ a root of}\, P(v,z_n),\,(z',\xi)\in
\tilde{\Omega}\}.\]
\begin{lem}
The set $\tilde{W}$ is open, definable and $\pi_{(n\times p, n\times
0)}( \tilde{W})=\tilde{\Omega}$.
\end{lem}
\underline{Proof.} It is clear that $\tilde{W}$ is definable.\\
If $(z',z_n)\in\tilde{\Omega}$ then
$(z',z_n,\sigma(z_n,\ldots,z_n))\in\tilde{W}$, hence
$\tilde{\Omega}\subset\pi_{(n\times p, n\times 0)}( \tilde{W})$, the
other inclusion is trivial. Let us prove that $\tilde{W}$ is open.\\
Let $(z',z_n, v)\in \tilde{W}$, since $ (z',z_n)\in \tilde{\Omega}$, there exists $ r_1 >0$, such that the ball:
\begin{equation}
\{ (y',y_n)\in \mathbb{C}^n\,/\,|z'-y'| <r_1 ,\,|z_n-y_n|<r_1\}\subset\tilde{\Omega} .
\end{equation}
- If $v=0$, we have then $(z',0)\in \tilde{\Omega}$,
there exists $ r_2 >0$, such that the ball:
\begin{equation}
\{ (y',y_n)\in \mathbb{C}^n\,/\,|z'-y'| <r_2 ,\,|y_n|
<r_2\}\subset\tilde{\Omega} .
\end{equation}
We put $r=inf\{r_1, r_2\}$, let $\rho >0$, $\rho < (\frac{r}{2})^p$, and consider the ball:
\[ B_{(z',z_n, 0)}:= \{ (y',y_n,u)\in \mathbb{C}^n\times \mathbb{C}^p\,/\, |z'-y'| <r ,\,|y_n-z_n|
<r, \mid u\mid <\rho\}.\]
We have $B_{(z',z_n, 0)}\subset\tilde{W}$. Indeed, if $  (y',y_n,u)\in B_{(z',z_n, 0)}$, we have $ (y',y_n)\in \tilde{\Omega}$ by 3.1.\\
Let $\eta$ be a root of the polynomial $ P(u,z_n) : = z_n^p + \sum_{ j=1}^pu_jz_n^{p-j}$, where $u= (u_1,\ldots,u_p)$.
By lemma 3.1, we have $ \mid \eta\mid\leq 2\,max_{j=1}^p|u_j|^{\frac{1}{j}}$. But, we have $ |u_j|\leq |u|\leq (\rho^{\frac{1}{p}})^j$, hence
$\mid \eta\mid \leq 2\rho^{\frac{1}{p}} < r$, which proves $ (y',\eta)\in \tilde{\Omega}$ by 3.2 and then $(y',y_n,u)\in \tilde{W}$ .\\

- If $v\neq 0$, we denote by $\xi_1,\ldots,\xi_p$ the roots of the polynomial $P(v,z_n) := z_n^p+\sum_{j=1}^pv_jz_n^{p-j}$, where $v=(v_1,\ldots,v_p)$.
Since $(z',z_n, v)\in \tilde{W}$, for each $j=1,\ldots,p$ there exists $\rho_j >0$ such that:
\begin{equation}
 \{ (y',y_n)\in \mathbb{C}^n\,/\, \mid y'-z'\mid < \rho_j,\,\mid z_n-\xi_j\mid < \rho_j\}\subset \tilde{\Omega} .
 \end{equation}
We put $r= min\{ r_1,\rho_1,\ldots, \rho_p\}$, where $r_1$ is the real in 3.1, and $\rho=max_{j=1}^pl_j$, where $l_j=(2\mid v\mid)^{\frac{1}{j}}$.\\
Let $\epsilon >$ such that $2\rho\epsilon <r$ and $\forall j=1,\ldots,p,\,\,\,\rho^j\epsilon^p\leq \mid v\mid$.
We  consider the ball;
\[ B_{ (z',z_n, v)}:= \{(y',y_n,u)\in \mathbb{C}^n\times \mathbb{C}^p\,/\, |z'-y'| <r ,\,|y_n-z_n|
<r, \mid u_j - v_j\mid <\rho^j\epsilon ^p,\,j=1\ldots,p\}.\]
We have $B_{ (z',z_n, v)}\subset \tilde{W}$. Indeed, if $(y',y_n,u)\in B_{ (z',z_n, v)}$, we see that $(y',y_n)\in \tilde {\Omega}$ by 3.1.\\
Let  $\eta$ be  a root of the polynomial $
P(u,z_n)$, where $u=(u_1,\ldots,u_p)$,  we have \[ \prod_{j=1}^p(\eta-\xi_j)= \eta^p +
\sum_{j=1}^pv_j\eta^{p-j}.\] Since $\eta^p + \sum_{j=1}^p
u_j\eta^{p-j}=0$, we have \[ \prod_{j=1}^p(\eta-\xi_j)= \sum_{j=1}^p
( v_j - u_j)\eta^{p-j} ,\] hence \[\prod_{j=1}^p|(\eta-\xi_j)|\leq
\sum_{j=1}^p
\mid v_j - u_j\mid \mid\eta\mid^{p-j}.\]
By lemma 3.1,$\mid \eta\mid \leq 2 max_{j=1}^p\mid u_j\mid$, since $\mid u_j\mid\leq \mid u_j-v_j\mid + \mid v_j\mid$, we deduce that $ \mid \eta\mid < 2\rho$. Hence:
\[ \prod_{j=1}^p|(\eta-\xi_j)|< 2^p\epsilon^p \rho^p .\]
There exists $\xi_k$ such that $\mid \xi_k-\eta\mid < 2\rho \epsilon < r$, hence $(y',\eta)\in \tilde {\Omega}$ by 3.3, which proves the inclusion.
\begin{rem}
Let $ \psi:\tilde{\Omega}\times \mathbb{C}^p\rightarrow
\tilde{\Omega}\times \mathbb{C}^p$ be the mapping defined by
$\psi(z,v)= (z,\sigma(v))$.
Clearly, $\psi$ is surjective and generically a submersion.\\
We put $\tilde{U}=\psi^{-1}(\tilde{W})\subset \mathbb{C}^n\times
\mathbb{C}^p$, we see that $\tilde{U}$ is open, symmetric and
$\psi(\tilde{U})= \tilde{W}$. We note also that for each
$(z',z_n,v_1,\ldots,v_p)\in\tilde{U}$ we have
$(z',v_j)\in\tilde{\Omega},\,1\leq j\leq p$.
\end{rem}

 In the following, we keep the notations and definitions above.
 \section{ Division with generic polynomial}
\begin{thm}
Let $P(v,z_n)$ be a generic polynomial in $z_n$ of degree $p$. If
$g\in\mathcal{A}(\Omega)$, then there exist an open definable set
$W\subset \Omega\times \mathbb{R}^p$, $ q\in\mathcal{A}(W)$ and
$h_j\in\mathcal{A}(\pi_{(n\times p,n-1\times p)}(W)),\,1\leq j\leq
p$, such that, for each $(x,t)\in W$,
\[ g(x)= q(x,t)P(t,x_n) +\sum_{j=1}^ph_j(x',t)x_n^{p-j} .\]
\end{thm}
\underline{Proof}. There is an open set $\tilde{\Omega}\subset\mathbb{C}^n$ and $\tilde{g}$ a holomorphic function on $\tilde{\Omega}$, such that,
$\tilde{\Omega}\cap\mathbb{R}^n = \Omega$ and $\tilde{g}_{\mid\Omega}=g $.
If $(z',z_n,v_1,\ldots,v_p)\in\tilde{U}$, we
have:
\[ \tilde{g}(z)-\tilde{g}(z',v_1)= (z_n-v_1)\tilde{g_1}(z,v),\]
with $\tilde{g_1}$ is a definable holomorphic function on
$\tilde{U}$. By repeating this process with $\tilde{g_1}$, we have:
\[ \tilde{g_1}(z,v)-\tilde{g_1}(z',v_2,v)= (z_n-v_2)\tilde{g_2}(z,v),\]
with $\tilde{g_2}$ is again a definable holomorphic function on
$\tilde{U}$. At the end, we get:
\[ \tilde{g}(z)=
\tilde{q}(z,v)(z_n-v_1)\ldots(z_n-v_p)+\sum_{j=1}^p\tilde{h_j}(z',v)z_n^{p-j},\]
with $\tilde{q}$ definable and holomorphic on $\tilde{U}$, and
$\tilde{h_j},\, 1\leq j\leq p$, definable and holomorphic on
$\pi_{(n\times p, n-1\times p)}(\tilde{U})$. The functions
$\tilde{q}$ and $\tilde{h_j},\, 1\leq j\leq p$ are symmetric by
construction. According to holomorphic Newton's theorem [
\cite{Toug}, IX, 6.5], there are holomorphic functions $\tilde{Q}$
on $\tilde{W}$  and $\tilde{H_j}$ on $\pi_{(n\times p, n-1\times
p)}(\tilde{W}),\, 1\leq j\leq p$, such that $\tilde{q}=
\tilde{Q}\circ\psi$ and $\tilde{h_j}=\tilde{H_j}\circ\psi,\, 1\leq
j\leq p$. Since $\psi$ is surjective, these functions are definable.
Since $P\circ \psi(z,v)=(z_n-v_1)\ldots(z_n-v_p)$ and $\tilde{g}\circ\psi(z,v)=\tilde{g}(z)$, we see that
\[ \tilde{g}(z)=
\tilde{Q}(z,v)P(v,z_n)+\sum_{j=1}^p\tilde{H_j}(z',v)z_n^{p-j} ,\]
since $\psi$ is
generically a submersion.
By restriction the above relation to $W:=\tilde{W}\cap \Omega\times
\mathbb{R}^p$ and taking the real parts, we get
\[ g(x)= Re\tilde {Q}(x,t) P(t,x_n) \sum_{j=1}^p Re \tilde{H_j}(x',t)x_n^{p-j},\,\,\forall (x,t)\in W,\]
which proves the theorem by putting $q(x,t)= Re\tilde {Q}(x,t)$ and
$h_j(x',t)= Re \tilde{H_j}(x',t),\,1\leq j\leq p.$

\section{ Semi-Global Weierstrass Division Theorem}
As we have seen for the algebra $\mathcal{O}(\Omega)$, we consider
the open set $\Omega$ equipped with the topology induced by
$SM\mathcal{A}(\Omega)$. This topological space will be denoted by
$\Omega(\mathcal{A})$. It is clear that $\Omega(\mathcal{A})$ is
Noetherian, i.e. every decreasing sequence of closed sets
stabilizes. If $A\subset \Omega$, we set:
\[ I(A)=\{ f\in\mathcal{A}(\Omega) /\,f(x)=0,\,\forall x\in A\}.\]
It is clear that $I(A)$ is an ideal.  Let
$ Z\subset \Omega(\mathcal{A})$ be a closed irreducible set, then $
I(Z)$ is a prime ideal. We put $\Lambda=
\frac{\mathcal{A}(\Omega)}{I(Z)}$, $\Lambda$ is a domain.
We denote by $ \Lambda[[x_1,\ldots,x_n]]$ the ring of formal series in the variables $x_1,\ldots,x_n$ with coefficients in  $\Lambda$. \\
Each $f\in\Lambda[[x_1,\ldots,x_n]]$  can be written as $
f=\sum_{n=0}^\infty f_n$ where $f_n\in\Lambda[x_1,\ldots,x_n]$ is a
homogenous polynomial of degree $n$.  We denote by  $\omega(f)$ the
smallest integer $n\in\mathbb{N}$, such that $f_n\neq 0$ and we call
$f_{\omega(f)}$ the  initial form of $f$. A formal series
$f\in\Lambda[[x_1,\ldots,x_n]]$ is called regular of order $p$ with
respect to $x_n$, if $ f(0,\ldots,0,x_n)= \delta x_n^p +
\ldots\ldots  $ , with $\delta\neq 0$ in $\Lambda =
\frac{\mathcal{A}(\Omega)}{I(Z)}$.
\begin{lem}
Let $g_1,\ldots,g_r\in\Lambda[[x_1,\ldots,x_n]]-\{0\}$, and put $p_j
=\omega(g_j),\,1\leq j\leq r$. Then there exists a linear
isomorphism $\sigma:\mathbb{R}^n\rightarrow \mathbb{R}^n$ such that
$g_j\circ\sigma$ is regular of order $p_j$ with respect to $x_n$,
for each $j= 1,\ldots r$.
\end{lem}
\underline{Proof}. Let $ P\in\Lambda[x_1,\ldots,x_n]$ be the product
of the initial forms of $g_1,\ldots,g_r$. $P$ is a homogenous
polynomial and  $P\neq 0$, hence there exists
$(\nu_1,\nu_2,\ldots,\nu_{n-1},1)\in\mathbb{Z}^n$ such that $
P(\nu_1,\nu_2,\ldots,\nu_{n-1},1)\neq 0$. We consider the linear
isomorphism \[\sigma(x_1,\ldots,x_n)= (x_1+\nu_1x_n,
x_2+\nu_1x_n,\ldots,x_{n-1}+\nu_{n-1}x_n, x_n).\] It is clear that
each $g_j\circ \sigma$ is regular of order $p_j$ with respect to
$x_n$, which proves the lemma.
\subsection{Semi-Global Weierstrass Preparation Theorem}
If $h\in \mathcal{A}(\Omega)$ we denote by $\overline{h}$ its image
by the canonical surjection   $ \mathcal{A}(\Omega)\rightarrow \Lambda = \frac{\mathcal{A}(\Omega)}{I(Z)}$. We
consider the morphism:
\[ \mathcal{A}(\Omega)\rightarrow \Lambda[[ x_1,\ldots,x_n]],\]
defined by
\[ g \rightarrow \sum_{\omega\in\mathbb{N}^n}\frac{\overline{D^\omega g}}{\omega !}x^\omega,\]
where $\omega=(\omega_1,\ldots,\omega_n),\,D^\omega g=\frac{
\partial^{|\omega|}}{ \partial x_1^{\omega_1}\ldots
\partial x_n^{\omega_n}}g,\,x^\omega= x_1^{\omega_1}\ldots x_n^{\omega_n}$,  and $\omega! =\omega_1 !\ldots\omega_n !$.\\

Let $\Omega\subset \mathbb{R}^n$ be a connected open definable set,
and $ g\in \mathcal{A}(\Omega)-\{0\}$. We consider $Z\subset
\Omega(\mathcal{A})$ a closed irreducible set. Put  $p\in\mathbb{N}$
the degree of the initial form of the series
$\sum_{\omega\in\mathbb{N}^n}\frac{\overline{D^\omega g}}{\omega
!}x^\omega$. By lemma 4.1, after a linear change of coordinates,
with coefficients in $\mathbb{Z}$, we can assume that $D^{
(0,\ldots,0,\nu)}g\in I(Z),\, 1\leq \nu < p$ and $\delta:= D^{
(0,\ldots,0,p)}g\notin I(Z)$. For each $a\in Z- V(\delta)$,
the germ of $g$ at $a$ is regular of order $p$ with respect to $x_n$, i.e. $g(a)=\frac{\partial}{\partial x_n}g(a)=\ldots = \frac{\partial^{p-1}}{\partial x_n^{p-1}}g(a)=0$ and $\frac{\partial^{p}}{\partial x_n^{p}}g(a)\neq 0$.\\
There exists  an open definable set $U\subset \Omega$,   such that $
Z- V(\delta)\subset U$ and $\delta(x)\neq 0,\,\forall x\in U$. By
theorem 3.1 , there exists an open definable set $W\subset
U\times\mathbb{R}^p$ with $\pi_{(n\times p,n\times 0)}( W)=U$, such
that:
\[ g(x)= q(x,t) P(t,x_n) \sum_{j=1}^p  h_j(x',t)x_n^{p-j},\,\,\forall (x,t)\in W,\eqno{(*)}\]
where $q\in\mathcal{A}(W),\,h_j\in \mathcal{A}(\pi_{(n\times p,n-1\times p)}(W)),\, 1\leq j\leq p$.\\
We consider the application $\varphi: W\subset U\times
\mathbb{R}^p\rightarrow  \mathbb{R}^n\times \mathbb{R}^p,$ defined
by \[\varphi(x,t)= (x,h_1(x',t),\ldots,h_p(x',t)).\] we put \[Y:=\{
(a',a_n,\sigma(a_n,\ldots,a_n))\in \mathbb{R}^n\times\mathbb{R}^p /\,(a',a_n)\in Z-V(\delta)\},\] by the
construction of $W$, we have $Y\subset W$.
\begin{rem}
Since $g$ is regular of order $p$ with respect to $x_n$ at any point
of $Y$, we deduce from $(*)$ that $\forall (a,b)\in Y,\, q(a,b)\neq
0$ and $h_j(a',b)=0,\,1\leq j\leq p$. We have then $\varphi( Y) = Z-
V(\delta)\times \{0\}$.
\end{rem}
\begin{lem}
The mapping $\varphi$ is a local diffeomorphism at every point of
$Y$. In addition, there exists an open neighborhood of $Y$,
$W'\subset W$, such that the restriction of $\varphi$ to $W'$ is
injective.
\end{lem}
\underline{Proof}.
 By shrinking the open set $W$, we can assume that $\forall (x,y)\in W$, $ q(x,y)\neq 0$. We note that
\[ \frac{\partial h_l}{\partial t_j}(a',b)=-\delta_{lj}q(a,b), \,\,\, j,l = 1,\ldots,p,(a',a_n,b)\in Y,\]
where $\delta_{lj}$ is the Kronecker symbol, hence
 the determinant of the matrix $[\frac{\partial h_l}{\partial t_j}(a',b)]_{j,l = 1,\ldots,p}$ is exactly $(-1)^p(q(a,b))^p$, we deduce the first assertion of the lemma.\\
 To prove the second assertion,
we consider the application  $\theta: W\subset U\times
\mathbb{R}^n\rightarrow \mathbb{R}^n\times \mathbb{R}^p$  defined by
\[\theta(x',x_n,t)= ( x',x_n, t_1+\sigma_1(x_n,\ldots,x_n),\ldots,
t_p +(-1)^{p-1}\sigma_p(x_n,\ldots,x_n)).\]
It is clear that $\theta (Y)= Z- V(\delta)\times\{0\}$ and  $\theta $ is a global diffeomorphism of W onto its image $ \theta(W)= W_1$. We consider its inverse $\theta^{-1}: W_1\rightarrow W$, and we put $\varphi_1= \varphi\circ\theta^{-1}$. We note that, for all $a\in Z-V(\delta)$, $\varphi_1(a,0)=(a,0)$ and $\pi_{(n\times p, n\times 0)}\circ\varphi_1 = \pi_{(n\times p, n\times 0)}$. It is clear  that to show the second assertion of the lemma, it suffices to show that there exists $W_1'\subset W_1 $ an  open neighborhood of $Z- V(\delta)\times\{0\}$ such that ${\varphi_1}_{|W'_1}$ is injective.\\ For each $a\in Z- V(\delta)$ there exists a ball $ B((a,0), r_a)\subset W_1$, such that the restriction of $\varphi_1$ to $ B((a,0), r_a)$ is injective. We set $W'_1=\bigcup_{a\in Z- V(\delta)}B((a,0), r_a)$, then the restriction of $\varphi_1$ to $W'_1$ is injective. Indeed, if $(x,y),\,(x',y')\in W'_1$ are such that $\varphi_1(x,y)= \varphi_1(x',y')$, then $x=x'$. Secondly there  exist $a , b \in Z- V(\delta)$ such that $ (x,y)\in B((a,0), r_a)$ and $ (x,y')\in B((b,0), r_b)$. If $r_a\leq r_b$, then $ (x,y)\in B((b,0), r_b)$ and since the restriction of $\varphi_1$ to the ball $ B((b,0), r_b)$ is injective, we deduce the result.\\

We put $W_1':=\varphi(W')$. The mapping $ \psi:\pi_{(n\times p,
n-1\times p)}(W')\rightarrow \pi_{(n\times p, n-1\times p)}(W_1')$
defined by
\[ \psi(x',t)= (x',h_1(x',t),\ldots,h_p(x',t))\]
is a global diffeomorphism with components in $\mathcal{A}(\pi_{(n\times p, n-1\times p)}(W'))$. We set $\gamma= \psi^{-1}$, we see then $\gamma=(x',\gamma_1,\ldots,\gamma_p)$, with $\gamma_j\in \mathcal{A}(\pi_{(n\times p, n-1\times p)}(W_1')),\,1\leq j\leq p$.\\
For each  $(x',0)\in\pi_{(n\times p, n-1\times p)}(W_1')$, we put $\gamma_j(x',0)=\psi_j(x'),\, 1\leq j\leq p$, and $\psi(x')=(\psi_1(x'),\ldots,\psi_p(x')).$ We see that $\forall (x',0)\in \pi_{(n\times p, n-1\times p)}(W_1'),\,h_j(x',\psi(x'))=0,\,1\leq j\leq p $.\\
By equation $(*)$, we have then in an open definable set   $
\Omega_1\supset Z-V(\delta)$:
\[  g(x) =q(x',\psi(x'))P(\psi(x'),x_n),\]
with $ \forall x\in \Omega_1,\,q(x',\psi(x'))\neq 0$.\\
If $A\subset\mathbb{R}^n$, then we denote the image of $A$ by the projection into the first $n-1$ coordinates by $ A_{n-1}$.\\
The following theorem summarizes the result we have just obtained.
\begin{thm}
Let $g\in\mathcal{A}(\Omega)$ and $ Z\subset \Omega(\mathcal{A})$ be
a closed irreducible set. Then there exist $p\in\mathbb{N}$,
$\delta\in\mathcal{A}(\Omega)-I(Z)$, and $ \Omega_1\subset\Omega$ an
open definable set such that:
\begin{enumerate}
\item $ Z- V(\delta)\subset \Omega_1$.
\item $g(x)= U(x)( x_n^p+\sum_{j=1}^p r_j(x')x_n^{p-j})$ where $U\in \mathcal{A}(\Omega_1),\, U(x)\neq 0,\,\forall x\in \Omega_1$,
$r_j\in \mathcal{A}(({\Omega_1})_{n-1}),\, r_j(a')=0,\,\forall a'\in (Z- V(\delta))_{n-1},\,1\leq j\leq p$.\\
We say that $g$ is equivalent on $ \Omega_1$ to the polynomial $
x_n^p+\sum_{j=1}^p r_j(x')x_n^{p-j}$.
\end{enumerate}
\end{thm}
\begin{cor}
Let $ Z\subset \Omega(\mathcal{A})$ be a closed irreducible set and
$g_1,\ldots,g_s\in \mathcal{A}(\Omega)$. Then there exist
$p_1,\ldots p_s\in\mathbb{N}$, $\delta\in \mathcal{A}(\Omega)-I(Z)$,
and $\Omega_1\subset \Omega$ an open definable set such that:
\begin{enumerate}
\item $ Z- V(\delta)\subset \Omega_1$.
\item $g_l(x)= U_l(x)( x_n^{p_l}+\sum_{j=1}^{p_l} r_{j,l}(x')x_n^{p_l-j})$ where $U_l\in \mathcal{A}(\Omega_1),\, U_l(x)\neq 0,\,\forall x\in \Omega_1$,
$r_{j,l}\in \mathcal{A}(({\Omega_1})_{n-1}),\,
r_{j,l}(a')=0,\,\forall a'\in (Z- V(\delta))_{n-1},\,1\leq j\leq
p_l$,  $l=1,\ldots,s$.
\end{enumerate}
\end{cor}
\underline{Proof}. We put $p_l,\,1\leq l\leq s$, the degree of the
initial form of the series
$\sum_{\omega\in\mathbb{N}^n}\frac{\overline{D^\omega g_l}}{\omega
!}x^\omega\in \frac{\mathcal{A}(\Omega)}{I(Z)}[[x_1,\ldots,x_n]]$.
By lemma 4.1, after a linear change of coordinates, we can suppose
that $D^{(0,\ldots,0,\nu)}g_l\in I(Z)$ if $\nu <p_l$
and $\delta_l:=D^{(0,\ldots,0,p_l)}g_l\notin I(Z),\,1\leq l\leq s.$\\
By theorem 4.4, for each $l=1,\ldots,s,$ there exists
$\Omega_l\subset\Omega$ an open definable such that $ Z-
V(\delta_l)\subset \Omega_l$ and $g_l(x)= U_l(x)(
x_n^{p_l}+\sum_{j=1}^{p_l} r_{j,l}(x')x_n^{p_l-j})$ where $U_l\in
\mathcal{A}(\Omega_l),\, U_l(x)\neq 0,\,\forall x\in \Omega_l$,
$r_{j,l}\in \mathcal{A}(({\Omega_l})_{n-1}),\, r_{j,l}(a')=0,\,\forall a'\in (Z- V(\delta_l))_{n-1},\,1\leq j\leq p_l$,  $l=1,\ldots,s$.\\
We put $\delta=\delta_1\ldots\delta_s$ and $\Omega_1=
\cap_{l=1}^s\Omega_l$. It is easy to see that $\delta$ and
$\Omega_1$, verify the corollary.
\begin{cor}
Let $g_1,\ldots,g_s\in\mathcal{A}(\Omega)$, then there exists a
finite number of definable open subsets of $\Omega$:
$\Omega_1,\ldots,\Omega_r$, such that,
\begin{enumerate}
\item $\Omega= \cup_{j=1}^r\Omega_j$,
\item $\forall l,\,1\leq l \leq s$, the restriction of $g_l,\,$ on $\Omega_j,\,1\leq j\leq r$, is equivalent to a unitary polynomial with respect to $x_n$ with coefficients in
$\mathcal{A}((\Omega_j)_{n-1})$.
\end{enumerate}
\end{cor}
\underline{Proof}. Recall that $\Omega(\mathcal{A})$ is a Noetherian
topological space, hence $ \Omega(\mathcal{A})= Z_1\cup
Z_2\cup\ldots\cup Z_\nu$, where $ Z_i$ is
a closed irreducible set, $ 1\leq i\leq \nu$. Applying the preceding corollary to $g_1,\ldots,g_s$ and $Z_1\subset \Omega$, there exist $\delta_1\in\mathcal{A}(\Omega)-I(Z_1)$ and $ \Omega_1\subset\Omega$ an open definable set such that:$ Z_1- V(\delta_1)\subset \Omega_1$ and the restriction of each  $g_l,\,1\leq l\leq s$, to $\Omega_1$ is equivalent to a unitary polynomial with respect to $x_n$ with coefficients in $\mathcal{A}((\Omega_j)_{n-1})$.\\
Now consider the decomposition of the closed  set $Z_1\cap
V(\delta_1)$ into irreducible sets: $Z_1\cap V(\delta_1)=
Z_{1,1}\cup Z_{1,2}\cup\ldots Z_{1, s_1}$. We repeat the same thing
with $Z_{1,\mu},\, 1\leq \mu\leq s_1$, and $g_1,\ldots,g_s$. This
process stops because the space is Noetherian. We do the same work
with $ Z_2,\ldots, Z_\nu$, which proves the corollary.

\section{$\mathcal{A}$-semi-analytic sets}
We now come to the setting of global semi-analytic sets. Let
$\Omega\subset \mathbb{R}^n$ be an open definable set.
\begin{defin}
A set $ A\subset \Omega$ is called $\mathcal{A}$- semi analytic set
if there exists a finite number of open definable sets
$\Omega_j\subset \Omega,\,1\leq j\leq
s,\,\Omega=\bigsqcup_{j=1}^s\Omega_j$ such that $ A\cap\Omega_j$ is
a finite union of sets of the form:
\[ \{ x\in\Omega /\, \varphi_0(x)=0,\,\varphi_1(x)>0,\ldots,\varphi_k >0\},\]
with $\varphi_0,\varphi_1,\ldots,\varphi_k\in\mathcal{A}(\Omega_j)$.
\end{defin}
Note that the $\mathcal{A}$-semi analytic sets form a boolean algebra of subsets of $\Omega$.\\
Before state and show the main result of this paper, we recall a
result of [ \cite{LVD},ch.2, theorem 2.7], see also [ \cite{Coste},
theorem 2.3.1] in the case of polynomials and [\cite{Bier} , theorem
2.6 ]
in the local analytic setting. First we recall some notations.\\
Let $X$ be a nonempty topological space and $\bf{E}$ a ring of continuous real-valued functions $f: X\rightarrow \mathbb{R}$.\\
Call a set $ A\subset X$ an $\bf{E}$-set if $A$ is finite union of
sets of the form:
\[\{x\in X /\,\varphi_0(x)=0,\varphi_1(x) >0,\ldots,\varphi_s(x) >0\},\,\,\mbox{ with}\,\,\varphi_0,\varphi_1,\ldots,\,\,\varphi_s\in \bf{E} .\]
Let $ A\subset X$, if $f,g :X\rightarrow \mathbb{R}$ are functions
such that $\forall x\in A,\,f(x) < g(x)$, we put $ (f,g): =
\{(x,t)\in A\times \mathbb{R}/\, f(x) < t < g(x)\}$.
\begin{thm} \cite{LVD}
Let $\varphi_1(T),\ldots,\varphi_M(T)\in {\bf{E}}[T]$. Then the list
$\varphi_1,\ldots,\varphi_M$ can be augmented to a list
$\varphi_1,\ldots,\varphi_N$ in ${\bf{E}}[T]$ ($M\leq N)$, and $X$
can be partitioned into finitely many $\bf{E}$-sets $ X_i,\,1\leq
i\leq k$, such that for each connected component $C$ of each $X_i$
there are continuous real-valued functions $\xi_{C,1} < \ldots <
\xi_{C,m_C}$ on $C$ with the following two properties:
\begin{enumerate}
\item each function $\varphi_n,\, 1\leq n\leq N$, has constant sign ($-1, 0, 1)$ on each of the graphs $\Gamma(\xi_{C,j}),\, 1\leq j\leq \mu(C)$, and on each of the sets
$(\xi_{\xi_{C,j}},\xi_{\xi_{C,j+1}}),\, 0\leq j\leq \mu(C)$, where
$\xi_{C,0}:=-\infty$ and $\xi_{C,j+1}:=\infty$ are constant
functions on $C$ by convention,
\item each of the sets $\Gamma(\xi_{C,j})$ and $(\xi_{\xi_{C,j}},\xi_{\xi_{C,j+1}})$ from $(1)$ is of the form
\[ \{(x,t)\in C\times\mathbb{R}/\, sign(\varphi_n(x,t)=\epsilon(n)\,\,\mbox{for}\,\,n = 1,\ldots,N \}\]
for a suitable sign condition $\epsilon:\{1,\ldots, N\}\rightarrow
\{ -1, 0, 1\}$.
\end{enumerate}
\end{thm}
\begin{thm}
If $ \Omega\subset\mathbb{R}^n$ is an open definable set and
$A\subset \Omega$ an $\mathcal{A}$- semianalytic set. Then $ A$ has
only a finitely many connected component, and each component is also
an $\mathcal{A}$- semianalytic set.
\end{thm}
\underline{Proof}. Since the structure $\mathcal{R}$ is o-minimal and $A$ is definable, we see that $A$ has only a finitely many connect components.\\
 For the second assertion of the theorem, we proceed  by induction on $n$. For $n=1$, the result is trivial, since $A$ is a finite union of
 points and intervals. We suppose $n >1$ and the result is true for $n-1$.\\
There exists a finite number of open definable sets $\Omega_j\subset
\Omega,\,1\leq j\leq s,\,\Omega=\bigsqcup_{j=1}^s\Omega_j$ such that
$ A\cap\Omega_j$ is a finite union of sets of the form:
\[ \{ x\in\Omega /\, \varphi_0(x)=0,\,\varphi_1(x)>0,\ldots,\varphi_k >0\},\]
with $\varphi_0,\varphi_1,\ldots,\varphi_k\in\mathcal{A}(\Omega_j)$.
We will restrict to each open $\Omega_j,\,1\leq j\leq s$. We can
then suppose that $A$ is described, in $\Omega$, by
$g_1,\ldots,g_s\in\mathcal{A}(\Omega)$. Applying the corollary 4.6
to this list, there exists a finite number of definable open subsets
of $\Omega$: $\Omega_1,\ldots,\Omega_r$, such that,
\begin{enumerate}
\item $\Omega= \cup_{j=1}^r\Omega_j$,
\item $\forall l,\,1\leq l \leq s$, the restriction of $g_l,\,$ on $\Omega_j,\,1\leq j\leq r$, is equivalent to a unitary polynomial with respect to $x_n$ with coefficients in
$\mathcal{A}((\Omega_j)_{n-1})$.
\end{enumerate}
For each $g_l,\,1\leq l\leq s$, we can then suppose that the
restriction of $g_l$ on each $\Omega_j,\, 1\leq j\leq r$, is in
$\mathcal{A}((\Omega_j)_{n-1})[x_n]$. The result follows from
theorem 5.2 by taking ${\bf{E}}=\mathcal{A}((\Omega_j)_{n-1})$   and
the induction hypothesis.
\section{Semi-global Weierstrass division theorem}
We give here a semi global Weierstrass division theorem for the
algebra $ \mathcal{A}(\Omega)$. We keep the notations and
assumptions of theorem 4.4. If
$\gamma,\gamma'\in\mathcal{A}(\Omega)$ and $\gamma$ is a multiple of
$\gamma'$ in $\mathcal{A}(\Omega)$, we write $ \gamma
\succsim\gamma'.$ Let $ Z\subset \Omega(\mathcal{A})$ be a closed
irreducible set and $\gamma\in\mathcal{A}(\Omega)- I(Z)$, we put
\[ \mathcal{A}_{Z,\gamma}=\{ f\,/\mbox{ There is an open definable set}\, W\subset\Omega,\, Z- V(\gamma)\subset W,\,\mbox{with}\,f\in \mathcal{A}(W)\}.\]
If $ \gamma \succsim\gamma'$, we have a restriction morphism $
\mathcal{A}_{Z,\gamma'}\rightarrow \mathcal{A}_{Z,\gamma}$. We put
\[ \mathcal{A}_{Z}= Lim_{\gamma\in\mathcal{A}(\Omega)-
I(Z)}\mathcal{A}_{Z,\gamma} .\]
\begin{thm}{[ Division theorem]}
Let $f\in \mathcal{A}(\Omega)$, then there exists an open definable
set $\Omega_1\subset \Omega$ such that $ Z- V(\delta)\subset
\Omega_1$ and $\forall x\in\Omega_1$:
\[ f(x)= q(x) g(x) +\sum_{j=1}^pr_j(x')x_n^{p-j} ,\]
where $q\in \mathcal{A}(\Omega_1)$, $r_j\in
\mathcal{A}((\Omega_1)_{n-1}),\,j\leq j\leq p$. The functions $ q,
r_j,\,j\leq j\leq p,$ are unique in $\mathcal{A}_Z$.
\end{thm}
\underline{Proof}. The uniqueness follows trivially from the fact
that  the germ of $g$ at each point $a\in Z-V(\delta)$ is regular of
order $p$ with respect to $x_n$.\\
By  theorem 4.4, there exists an open set $\Omega_1\subset \Omega$,
such that $ Z- V(\delta)\subset \Omega_1$ and $g$ is equivalent, in
$\mathcal{A}(\Omega_1)$,  to the polynomial  $x_n^p
+\sum_{j=1}^pr_j(x')x_n^{p-j}$. By theorem 3.1, we can divide $f$ by
the generic polynomial $x_n^p+\sum_{j=1}^pt_jx_n^{p-j}$, and we
substitute   $(r_1(x'),\ldots,r_p(x'))$ in the place of $
(t_1,\ldots,t_p)$, which gives the result.
\subsection*{Acknowledgment}
the results of this paper were completed during my stay in the ICTP centre . I want to thank the mathematics section for the invitation.

\end{document}